\documentstyle[12pt]{article}
\author{David Ruelle and  Amie Wilkinson}
\title{Absolutely singular dynamical foliations}
\newtheorem{theorem}{Theorem}[section]

\newtheorem{cor}[theorem]{Corollary}

\def\proof{{\bf {\medskip}{\noindent}Proof: }}

\def\remark{{\bf {\bigskip}{\noindent}Remark: }}

\def\eproof{$\Box$ \bigskip}

\def\title{\em}

\def\mod{\hbox{mod }}

\def\dim{\hbox{dim}}
\def\exp{\hbox{exp}}

\def\Diff{\hbox{Diff}}

 \newlength{\figboxwidth} 

\begin{document}
\maketitle

\section*{Introduction}

Let $A_2$ be the automorphism of the $2$-torus,
${\bf T}^2={\bf R}^2\slash{\bf Z}^2,$
given by $\left(\begin{array}{cc}
2 & 1 \\
1 & 1
\end{array}
\right).$
Let $A_3$ be the automorphism of 
the $3$-torus ${\bf T}^3={\bf R}^3\slash{\bf Z}^3$ 
given by $\left(\begin{array}{cc}
A_2 & 0\\
0 & 1
\end{array}
\right).$
Let $\Diff^2_\mu({\bf T}^3)$ be the set of $C^2$
diffeomorphisms of ${\bf T}^3$ that
preserve Lebesgue-Haar measure $\mu$.

In \cite{SW}, M. Shub and A. Wilkinson prove the following theorem.

\bigskip

\noindent{\bf Theorem: }{\em Arbitrarily close to $A_3$ there is a $C^1$-open set
$U \subset \Diff^2_\mu({\bf T}^3)$ such
that for each $g\in U$,
\begin{enumerate}
\item $g$ is ergodic. 
\item There is an equivariant fibration
$\pi:{\bf T}^3\to {\bf T}^2$  such that $\pi g=A_2 \pi$
The fibers of $\pi$ are the leaves of a foliation
${\cal{W}}^c_g$ of ${\bf T}^3$ by $C^2$ circles. In particular,
the set of periodic leaves is dense in ${\bf T}^3$.
\item There exists $\lambda^c>0$ such that,
for $\mu$-almost every $w\in {\bf T}^3$,  if $v\in T_w{\bf T}^3$
is tangent to the leaf of ${\cal{W}}^c_g$ containing $w$, then
$$\lim_{n\to \infty} {1\over n}\log \| T_wg^n v\| =\lambda^c.$$
\item Consequently, there exists a set $S\subseteq {\bf T}^3$
of full $\mu$-measure that meets every
leaf of ${\cal{W}}^c_g$ in a set of leaf-measure $0$. The 
foliation ${\cal{W}}^c_g$ is
not absolutely continuous.
\end{enumerate}}

\medskip

\noindent Additionally, it is shown that the diffeomorphisms
in $U$ are nonuniformly hyperbolic and Bernoullian.
In this note, we prove:

\medskip
\noindent{\bf Theorem I: } {\em Let $g$ satisfy conclusions
1.--3. of the previous theorem.
Then  there exist  $S\subseteq {\bf T}^3$ of full $\mu$-measure 
and $k\in{\bf N}$ such that 
$S$ meets every
leaf of ${\cal{W}}^c_g$ in exactly $k$ points.
The foliation ${\cal{W}}^c_g$ is
{\em absolutely singular.}}

\medskip

\remark In A. Katok's example of an
absolutely singular foliation in \cite{M},
the leaves of the foliation meet the set of full measure in
one point.  In the \cite{SW} examples, the set 
$S$ may necessarily meet leaves of ${\cal{W}}^c_g$ in more
than one point, as the following argument of Katok's shows.

It follows from Theorem II in \cite{SW2} that for $k\in{\bf Z}_+$ and
for small $a,b>0$, the
map $g = j_{a,k}\circ h_{b}$ satisfies the hypotheses of Theorem I, where 
$$h_b(x,y,z)= (2x+y, x+y, x+y+z + b\sin 2\pi y),\qquad\hbox{and}$$
$$j_{a,k}(x,y,z) = (x,y,z) + a \cos (2\pi k z)\cdot (1+\sqrt{5}, 2, 0).$$

For $k\in{\bf N}$,
let $\rho_k$ be the vertical translation that sends $(x,y,z)$
to $(x,y,z+{1\over k})$.  Note that $h_{b}\circ \rho_k = \rho_k \circ h_b$
and $j_{a,k}\circ \rho_k = \rho_k\circ j_{a,k}$.  Thus  
$g\circ \rho_k = \rho_k\circ g$.

The fibration $\pi:{\bf T}^3\to{\bf T}^2$ was obtained in \cite{SW}
by using the persistence of normally hyperbolic submanifolds under 
perturbations.  In the present case the symmetries $\rho_k$ preserve the 
fibers of the trivial fibration $P:{\bf T}^3\to{\bf T}^2$ from which one 
starts, and also the maps $g$.  Therefore the fibers of 
$\pi:{\bf T}^3\to{\bf T}^2$ (i.e., the leaves of center foliation 
${\cal W}_g^c$) are invariant under the action of the finite group 
$<\rho_k>$.

Let $S$ be the (full measure) 
set of points in ${\bf T}^3$ for which the center direction
is a positive Lyapunov direction (i.e. for which
conclusion 3 holds).
Since $\rho_k({\cal{W}}^c_g) = {\cal{W}}^c_g$,
it follows that $\rho_k S=S$.
If $p\in S\cap {\cal W}^c(p)$, then  $\rho_k(p)\in \rho_k(S)\cap 
\rho_k({\cal{W}}^c(p)) = S\cap {\cal W}^c(p)$; that is, 
$S\cap {\cal W}^c(p)$ contains at least $k$ points.   

Thus Theorem I is ``sharp'' in the sense that we cannot say
more about the value of $k$ in general.  We see no
reason why $k=1$ should hold even for a residual set in $U$.

\bigskip

Theorem I has an interesting interpretation.
Recall that a {\em $G$-extension} of a dynamical system
$f:X\to X$ is a map  $f_\varphi:X\times G\to X\times G$, 
where $G$ is a compact group,
of the form $(x,y)\mapsto (g(x), \varphi(x)y)$.
If $f$ preserves $\nu$, and $\varphi:X\to G$ is measurable,
then $f_\varphi$ preserves the product of $\nu$ with
Lebesgue-Haar measure on $G$.  A ${\bf Z}\slash k{\bf Z}$-extension 
is also called a $k$-point extension.

	Let $\lambda$ be an invariant probability measure for a $k$-point 
extension of $f:X\to X$, and $\{\lambda_x\}$ the family of conditional 
measures associated with the partition $\{\{x\}\times G\}$.  We remark that 
if $\lambda$ is ergodic, then 
each atom of $\lambda_x$ must have the same weight 
$1/k$ (up to a set of $\lambda$-measure $0$).

	Now take $g\in U$.  Choose a coherent orientation on the leaves of 
$\{\pi^{-1}(x)\}_{x\in T^2}$.  Take $h:{\bf T}^3\to{\bf T}^2\times{\bf T}$ 
to be any continuous change of coordinates such that $h$ restricted to 
$\pi^{-1}(x)$ is smooth and orientation preserving to $\{x\}\times{\bf T}$.
We may then write $F=h\circ g\circ h^{-1}:{\bf T}^2\times{\bf T}\to
{\bf T}^2\times{\bf T}$ in the form
$$	F(x,p)=(A_2x,\varphi_x(p))      $$
where $\varphi_x:{\bf T}\to{\bf T}$ is smooth and orientation preserving.  If 
$P:{\bf T}^2\times{\bf T}\to{\bf T}^2$ is the projection on the first 
factor of the product, we have $P\circ h=\pi$.  Therefore, writing $\lambda
=h^*\mu$, we have $P^*\lambda=\pi^*\mu$.  Let $\{\lambda_x\}$ be the 
disintegration of the measure $\lambda$ along the fibers $\{x\}\times
{\bf T}$.  By a further measurable change of coordinates, smooth along 
each $\{x\}\times{\bf T}$ fiber, we may assume that $\lambda$-almost 
everywhere, the atoms of $\lambda_x$ are at $l/k$, for $l=0,\ldots,k-1$.  
But then $\varphi_x$ permutes the atoms cyclically, and we obtain the 
following corollary.

\medskip
	
\noindent{\bf Corollary:}  {\sl For every $g\in U$ there exists $k\in
{\bf N}$ such that $({\bf T}^3,\mu,g)$ is isomorphic to an (ergodic)
$k$-point extension of $({\bf T}^2,\pi^*\mu,A_2)$.}
\bigskip

M. Shub has observed that if $g=j_{a,k}\circ h_b$, then 
$\pi^*\mu$ is actually Lebesgue measure on ${\bf T}^2$.

\section{Proof of Theorem I}

The proof of Theorem I follows from a more general result about
fibered diffeomorphisms.  Before stating this result,
we describe the underlying setup and assumptions.

Let $X$ be a compact metric space
with Borel probability measure $\nu$, and let $f:X\to X$ be
invertible and ergodic with respect to $\nu$.
Let $M$ be a closed Riemannian manifold
and  $\varphi:X\to\Diff^{1+\alpha} (M)$
a measurable map.  Consider the skew-product transformation
$F:X\times M\to X\times M$ given by
$$F(x,p) = (f(x), \varphi_x(p)).$$
Assume further that there is an $F$-invariant ergodic probability
measure $\mu$ on $X\times M$ such that $\pi_\ast \mu = \nu$, where
$\pi:X\times M\to X$ is the projection onto the first factor.

For $x\in X$, let $\varphi^{(0)}_x$ be the identity map on $M$
and for $k\in {\bf Z}$, define  $\varphi^{(k)}_x$ 
by $$\varphi^{(k+1)}_x  = \varphi_{f^k(x)}\circ \varphi^{(k)}_x.$$
Since the tangent bundle to $M$ is measurably trivial,
the derivative map of $\varphi$ along the $M$ direction
gives a cocycle $D\varphi: X\times M\times{\bf Z}\to GL(n,{\bf R})$, where
$n=\dim(M)$:
$$(x,p,k)\mapsto D_p\varphi^{(k)}_x.$$

Assume that $\log^+\|D\varphi\|_{\alpha}\in L^1(X\times M,\mu)$,
where $\|\,\cdot\,\|_\alpha$ is the $\alpha$-H{\"o}lder norm.
Let $\lambda_1<\lambda_2\cdots<\lambda_l$ be
the Lyapunov exponents of this cocycle; they exist
for $\mu$-a.e. $(x,p)$ by Oseledec's Theorem and
are constant by ergodicity.  We call these the {\em fiberwise exponents}
of $F$. Under the assumptions just described, we have
the following result.

\bigskip

\noindent{\bf Theorem II: }{\em  Suppose that $\lambda_l<0$. Then there exists
a set $S\subseteq X\times M$ and an integer $k\geq 1$ such that
\begin{itemize}
\item $\mu(S) = 1$
\item For every $(x,p)\in S$, we have $\#(S\cap \{x\} \times M) = k$.
\end{itemize}}

\medskip

This has the immediate corollary:

\bigskip

\noindent{\bf Corollary: }{\em Let $f\in \Diff^{1+\alpha}(M)$.  If $\mu$ is an ergodic measure with
all of its exponents negative, then it is concentrated on the orbit of a 
periodic sink.}

\medskip

The corollary has a simple proof using regular neighborhoods. 
Our proof
is a fibered version.  Theorem I is also a corollary of Theorem II.
For this, the argument is actually applied to the inverse of $g$,
which has negative fiberwise exponents, rather than to $g$ itself,
whose fiberwise exponents are positive. As we described
in the previous remarks, there is a measurable change
of coordinates, smooth along the leaves of
${\cal{W}}^c_g$ in which  $g^{-1}$ is expressed as
a skew product of ${\bf T}^2\times{\bf T}$.

\remark  Without the assumption that $f$ is invertible, Theorem II
is false. An example is described by Y. Kifer \cite{Ki}, which
we recall here.  Let $f:{\bf T}\to {\bf T}$ be a $C^{1+\alpha}$
diffeomorphism
with exactly two fixed points, one attracting and
one repelling. Consider the following random
diffeomophism  of ${\bf T}$: with probability $p\in(0,1)$, apply
$f$, and with probability $1-p$, rotate by an angle chosen
randomly from the interval $[-\epsilon, \epsilon]$. 

Let $X=(\{0,1\}\times{\bf T})^{\bf N}$.  To generate
a sequence of diffeomorphisms
$f_0,f_1,\ldots$ according to the above rule,
we first define
$\varphi:X\to \Diff^{1+\alpha}({\bf T})$ by
$$\varphi(\omega) =\cases{f & if $\omega(0) = (0,\theta)$,\cr
R_\theta & if $\omega(0) = (1,\theta)$,}$$
where $R_\theta$ is rotation through angle $\theta$.
Next, we let $\nu_\epsilon$
be the product of $p, 1-p$-measure on $\{0,1\}$ with
the measure on ${\bf T}$ that is uniformly distributed
on $[-\epsilon, \epsilon]$. Then corresponding to 
$\nu_\epsilon^{\bf N}$-almost
every element $\omega\in X$ is the sequence 
$\{f_k  = \varphi(\sigma^k(\omega))\}_{k=0}^\infty,$
where  $\sigma:X\to X$ is the one-sided shift $\sigma(\omega)(n) =
\omega(n+1)$.

Put another way, the random
diffeomorphism is generated by the (noninvertible) skew product
$\tau: X\times{\bf T}\to X\times{\bf T}$,
where $\tau(\omega,x) = (\sigma(\omega),
\varphi(\omega)(x))$.
An ergodic $\nu_\epsilon$-stationary measure for this random diffeomorphism
is a measure  $\mu_\epsilon$ on ${\bf T}$ such that
$\mu_\epsilon\times \nu_\epsilon^{\bf N}$ is $\tau$-invariant and ergodic.
Such measures always exist (\cite{Ki}, Lemma I.2.2), but,
for this example, there is an ergodic stationary
measure with additional special properties.

Specifically, for every $\epsilon>0$,
there exists an ergodic $\nu_\epsilon$-stationary 
measure $\mu_\epsilon$ on ${\bf T}$
such that, as $\epsilon\to 0$, 
$\mu_\epsilon \to \delta_{x_0}$, in the weak topology,
where $\delta_{x_0}$ is Dirac measure concentrated
on the sink $x_0$  for $f$.
From this, it follows that, as $\epsilon\to 0$, the fiberwise Lyapunov 
exponent for $\mu_\epsilon$ approaches $\log |f'(x_0)| < 0$, which is
the Lyapunov exponent of $\delta_{x_0}$.
Thus, for $\epsilon$ sufficiently small, the fiberwise 
exponent for $\tau$ with respect to $\mu_\epsilon$ is negative.
Nonetheless, it is easy to see that
$\mu_\epsilon$ for $\epsilon>0$ cannot be uniformly distributed on $k$ atoms;
if $\mu_\epsilon$ were atomic,
then $\tau$-invariance of $\mu_\epsilon\times \nu_\epsilon^{\bf N}$
would imply that, for every $x\in {\bf T}$, 
\begin{eqnarray*}
\mu_\epsilon(\{x\}) & = &
p\mu_\epsilon(\{f^{-1}(x)\})+(1-p)\int_{-\epsilon}^{\epsilon}
\mu_\epsilon(\{R_\theta(x)\}) d\theta\\
& = & p\mu_\epsilon(\{f^{-1}(x)\}),
\end{eqnarray*}
which is impossible if  $\mu_\epsilon$ has finitely many atoms.
In fact, $\mu_\epsilon$ can be shown to be absolutely
continuous with respect to Lebesgue measure (see \cite{Ki},
p. 173ff  and the references cited therein). Hence invertibility
is essential, and we indicate in the proof of Theorem II
where it is used.

\bigskip

\noindent{\bf Proof of Theorem II:}
We first establish the existence of fiberwise ``stable manifolds''
for the skew product $F$.  A general theory of stable manifolds for
random dynamical systems is worked out in (\cite{Ki}, Theorem V.1.6;
see also \cite{LQ});
since we are assuming that all of the
fiberwise exponents for $F$ are negative, we are faced with the simpler
task of constructing fiberwise regular neighborhoods for $F$ 
(see the Appendix by Katok and Mendoza
in \cite{KH}).  We outline a proof, following closely \cite{KH}.

\begin{theorem} (Existence of Regular Neighborhoods)
There exists a set $\Lambda_0\subseteq X\times M$ of full measure
such that for $\epsilon > 0$:
\begin{itemize}
\item There exists a measurable function 
$r:\Lambda_0\to (0,1]$ and a collection
of embeddings $\Psi_{(x,p)}:B(0,q(x,p))\to M$ such that $\Psi_{(x,p)}(0) = p$
and $\exp(-\epsilon)< r(F(x,p))/r(x,p) < \exp(\epsilon)$.
\item If $\varphi_{(x,p)} =  \Psi^{-1}_{F(x,p)}\circ \varphi_x \circ \Psi_{(x,p)}:
B(0,r(x,p))\to {\bf R}^n$, then $D_0\varphi_{(x,p)}$
satisfies 
$$\exp(\lambda_1 - \epsilon) \leq \| D_0\varphi_{(x,p)}^{-1}\|^{-1},
\|D_0\varphi_{(x,p)}\|\leq \exp(\lambda_l + \epsilon).$$
\item The $C^1$ distance $d_{C^1}(\varphi_{(x,p)}, D_0\varphi_{(x,p)})<\epsilon$
in $B(0,r(x,p))$.
\item There exist a constant $K>0$ and a measurable function 
$A:\Lambda_0\to {\bf R}$
such that for $y,z\in B(0,r(x,p))$,
$$K^{-1}d(\Psi_{(x,p)}(y), \Psi_{(x,p)}(z))\leq \|y-z\| \leq
A(x) d(\Psi_{(x,p)}(y),\Psi_{(x,p)}(z)),$$
with $\exp(-\epsilon) < A(F(x,p))/A(x,p) < \exp(\epsilon)$.
\end{itemize}
\end{theorem}

\proof See the proof of Theorem S.3.1 in \cite{KH}.  \eproof

\bigskip

Decompose $\mu$ into a system of fiberwise
measures $d\mu(x,p) = d\mu_x(p) d\nu(x)$.
Invariance of $\mu$ with respect to $F$ implies that, for
$\nu$-a.e. $x\in X$,
$${\varphi_x}_\ast \mu_x = \mu_{f(x)}.$$

\begin{cor}\label{C=manif} There exists a set 
$\Lambda\subseteq X\times M$, and real numbers
$R>0$, $C>0$, and $c<1$ such that
\begin{itemize}
\item[(1)] $\mu(\Lambda)>.5$, and, if $(x,p)\in \Lambda$, then $\mu_x(\Lambda_x)>.5$,
where $\Lambda_x = \{p\in M\,\vert\, (x,p)\in\Lambda\}$, 
\item[(2)] If $(x,p)\in \Lambda$ and $d_M(p,q)\leq R$, then
$$d_M(\varphi^{(m)}_x(p), \varphi^{(m)}_x(q)) \leq C c^m d_M(p,q),$$
for all $m\geq 0$.
\end{itemize} 
\end{cor}

\proof This follows in a standard way
from the Mean Value Theorem and Lusin's Theorem.\eproof

\bigskip

To prove Theorem II, it suffices to show that there is
a positive $\nu$-measure set $B\subseteq X$, such that for $x\in B$,
the measure $\mu_x$ has an atom, as the following argument shows. For $x\in X$, 
let $d(x) = \sup_{p\in M} \mu_x(p)$. Clearly $d$ is measurable,
$f$-invariant, and positive on $B$. Ergodicity of $f$ implies that
$d(x)=d>0$ is positive and constant for almost all $x\in X$.  
Let $S=\{(x,p)\in X\times M\,\vert\,
\mu_x(p)\geq d\}$.  Observe that $S$ is $F$-invariant, has measure
at least $d$, and hence has measure $1$.  
The conclusions of Theorem~II
follow immediately.  

Let $\Lambda$, $R>0$, $C>0$, and $c<1$ be given by
Corollary~\ref{C=manif}, and let  $B=\pi(\Lambda)$.  Let $N$ be the number of 
$R/10$-balls needed to cover $M$.  
We now show that for  $\nu$-almost every $x\in B$, the measure $\mu_x$ has at least
one atom. 

For $x\in X$, let
$$ m(x) = \inf \sum  \hbox{diam }(U_j),$$
where the infimum is taken over all collections of closed balls
$U_1,\ldots, U_k$ in $M$ such that $k\leq N$ and  $\mu_x(\bigcup_{j=1}^k U_j) \geq .5$.
Let $m = \hbox{ess sup }_{x\in B} m(x)$. 

We now show that $m=0$.  If $m>0$, then there exists an integer $J$
such that
\begin{eqnarray}\label{e=finalsmall} 
C\Delta c^J N & < & m/2, 
\end{eqnarray}
where $\Delta$ is the diameter of $M$.  Let
${\cal{U}}$ be a cover of $M$ by $N$ closed balls of radius 
$R/10$.  For $x\in B$, let $U_1(x),\ldots,U_{k(x)}(x)$  be those
balls in ${\cal{U}}$ that meet $\Lambda_x$.  Since these balls
cover $\Lambda_x$, and $\mu_x(\Lambda_x)>.5$, it follows that
$\mu_x(\bigcup_{j=1}^{k(x)} U_j(x)) \geq .5$.  But 
${\varphi_x^{(i)}}_\ast \mu_x = \mu_{f^i(x)},$  and so it's
also true that
\begin{eqnarray}\label{e=largemeas}
\mu_{f^i(x)}(\bigcup_{j=1}^{k(x)} \varphi_x^{(i)}(U_j(x))) &\geq &.5,
\end{eqnarray}
for all $i$.

We now use the fact that $\varphi^{(i)}_x$ contracts regular
neighborhoods to derive a contradiction.  The balls
$U_j(x)$ meet $\Lambda_x$ and have diameter less than $R/10$,
and so by Corollary~\ref{C=manif}, (2), we have
\begin{eqnarray}\label{e=smalldiam}
\hbox{diam }({\varphi_x^{(i)}}(U_j(x))) &\leq & C \Delta c^i.
\end{eqnarray}
Let
$\tau:B\to {\bf N}$ be the first-return time of  $f^J$ to $B$,
so that $f^{J\tau(x)}(x)\in B$, and $f^{Ji}(x)\notin B$, for
$i\in \{1,\ldots,\tau(x)-1\}$.  Decompose the set $B$
according to these first return times:
$$B=\bigcup_{i=1}^\infty B_i\quad(\mod 0),$$
where $B_i=\tau^{-1}(i)$.  {\em Because $f$ is invertible and
$f^{-1}$ preserves measure}, we also have the mod 0 
equivalence:
$$B' := \bigcup_{i=1}^{\infty} f^{Ji}(B_i) = B\quad(\mod 0).$$

Let $y\in B'$.  Then $y=f^{Ji}(x)$, where $x\in B_i\subseteq B$, for
some $i\geq 1$.  It follows from the definition of $m(y)$ 
and inequalities (\ref{e=largemeas}), (\ref{e=smalldiam})
and (\ref{e=finalsmall}) that 
\begin{eqnarray*}
m(y) & \leq & \sum_{j=1}^{k(x)}\hbox{diam }(\varphi_x^{(Ji)}(U_j(x)))\\
&\leq & C  k(x) \Delta c^{Ji}\\
&\leq & C N \Delta c^{J}\\
& < & m/2.
\end{eqnarray*}
But then
\begin{eqnarray*}
m&=&\hbox{ess sup }_{x\in B} m(x)\\
& = & \hbox{ess sup }_{y\in B'} m(y)\\
& < & m/2,
\end{eqnarray*}
contradicting the assumption $m>0$.

Thus $m=0$, and, for $\nu$-almost every $x\in B$, we have
$m(x)=0$. If $m(x)=0$, then there is a sequence
of closed balls $U^1(x), U^2(x), \cdots$
with $\lim_{i\to\infty} \hbox{diam }(U^i(x))=0$ and
$\mu_x(U^i(x))\geq .5/N$, for all $i$.  Take
$p_i\in U^i(x)$; any accumulation point of $\{p_i\}$ is an atom for $\mu_x$.
Since we have shown that $\mu_x$ has an atom, for
$\nu$-a.e. $x\in B$,  the proof of Theorem~II is complete.\eproof

We thank Michael Shub and Anatole Katok for useful conversations.

\end{document}